\documentclass[review]{elsarticle}

\usepackage{amssymb,amsmath,amsthm}
\usepackage{color}
\usepackage{graphicx}
\usepackage{subfigure}
\usepackage{epstopdf}
\usepackage{graphics,epsfig}
\usepackage{algorithm}

\newcommand{\argmin}{\arg\!\min}
\renewcommand{\qedsymbol}{$\blacksquare$}

\journal{}

\begin{document}

\begin{frontmatter}

\title{Nonlinear Consensus Strategies for Multi-Agent Networks in Presence of Communication Delays and Switching Topologies: \\ Real-Time Receding Horizon Approach}

\author[mymainaddress]{Fei Sun\corref{mycorrespondingauthor}}
\author[mymainaddress]{Kamran Turkoglu\corref{mycorrespondingauthor}}

%
\cortext[mycorrespondingauthor]{Email: fei.sun@sjsu.edu, kamran.turkoglu@sjsu.edu}
%
\address[mymainaddress]{Aerospace Engineering, San Jos\'e State University, San Jose, CA 95192}

\begin{abstract}
This paper presents a novel framework which combines a non-iterative solution of Real-Time Nonlinear Receding Horizon Control (NRHC) methodology to achieve consensus within complex network topologies with existing time-delays and in presence of switching topologies. In this formulation, we solve the distributed nonlinear optimization problem for multi-agent network systems directly, \emph{in real-time}, without any dependency on iterative processes, where the stability and convergence guarantees are provided for the solution. Three benchmark examples on non-linear chaotic systems provide validated results which demonstrate the significant outcomes of such methodology.  
\end{abstract}

\begin{keyword}
multi-agent consensus problems\sep  nonlinear receding horizon control \sep
real-time optimization\sep switching topologies
\end{keyword}

\end{frontmatter}


\section{INTRODUCTION}

The important problem of finding distributed control laws governing network and/or corresponding multi-agent dynamics has been of interest within the control research community the for the recent decade. One important aspect of this research effort is to answer the existence and uniqueness properties of a distributed control law (and solution) that can lead each individual agent within the network to a common state value, which constitutes the core of the (multi-agent) network consensus problem. 

For several years, the complex and sophisticated nature of multi-agent consensus problem has been serving as a fruitful research field for the application of advanced control algorithms (such as formation flights, power grid networks, flock dynamics, cooperative control ... etc.) (\cite{saber2004}-\cite{zhao2015}). 

So far, in existing literature, consensus methodologies have been widely explored for multi-agent \emph{linear} dynamical systems in many important and ground-breaking studies (\cite{LiDuanChenHuang2010} - \cite{MovricLewis2014}). However, there exists some important studies which also investigate nonlinear consensus problems in multi agent systems. For example, Bausso et al. (2006) investigated a nonlinear protocol design that permits consensus on a general set of values \cite{BaussoGiarrePesenti2006}. In Qu et al. (2007) emphasis was put on nonlinear cooperative control for consensus of nonlinear and heterogeneous systems \cite{QuChunyuWang2007}. Liu et al. (2013) demonstrates the applicability of a variable transformation method to convert a general nonlinear consensus problem to a partial stability problem \cite{LiuXieRenWang2013}. Franco et al. (2008) investigates an important connection between nonlinear consensus problems and receding horizon methodology with existing time-delay information and constraints \cite{Franco2008}.  

In recent years, with the fascinating results improvements obtained in microprocessor related technologies, the applicability of real-time based control methodologies, and more specifically Receding Horizon Control based methodologies became more of a reality. In that sense, receding horizon control emerged as one of the existing control methodologies that could be adapted to consensus problem of multi-agent dynamics, for real-time solutions. Based on this approach, there have been many results developed for consensus problems and its applications. The work of \cite{dunbar2005} presented a distributed receding horizon control law for dynamically coupled nonlinear systems based on its linearization representative. The robust distributed receding horizon control methods were studied in \cite{LiShi2013} for nonlinear systems with coupled constraints and communication delays. \cite{shi2014} proposed a robust distributed model predictive control methods for nonlinear systems subject to external disturbances. In a very recent study Qiu and Duan (2014) investigated brain-storm type of optimization in combination with receding horizon control strategies for UAV formation flight dynamics, \cite{QiuDuan2014}.
 
One very interesting (and at the same time natural) extension to the problem is the consensus of nonlinear multi-agent networks with switching topologies and time delays. Li and Qu (2014) \cite{LiQu2014} provide an interesting approach to the finite time consensus problem of distributed nonlinear systems under general setting of directed and switching topologies. Jia and Tang (2012) \cite{JiaTang2012} worked on a specific directed graph topology to achieve consensus within nonlinear agents with switching topology and communication delays. In their studies, Ding and Guo (2015) \cite{DingGuo2015} concentrated their effort on sampled-data based leader-following consensus for nonlinear multi-agent systems with Markovian switching topologies and communication delays. 

In all of those existing valuable studies, one common fact remains as the iterative nature of the problem solution methodology, which enforces the iterative nature of the solution to still exist. This heavily deters the functionality of real-time solution methodologies.  
 
In this paper we aim to address this issue, and for this purpose we propose a novel framework which combines non-iterative solution of Real-Time Nonlinear Receding Horizon Control (NRHC) methodology to achieve consensus within complex network topologies with existing time-delays and in presence of switching topologies. With this we propose the following three novelties: (i) non-iterative solution is achieved in real-time which eliminates any need for an iterative algorithm, (ii) multi-agent consensus problem is solved in real-time for switching topologies and existing time-delays, (iii) stability and convergence of the consensus is guaranteed in existence of complex switching topologies, (iv) stability and convergence of the consensus is guaranteed in existence of inherent communication delays. 

The paper is organized as follows: In Section-\ref{sec:prob_formulation}, problem statement is defined. Distributed NRHC protocol is presented in Section-\ref{sec:dist_control} . Section-\ref{sec:stab_analysis} provides important results on convergence and stability assessment of the presented methodology, where in Section-\ref{sec:time-delay} multi-agent network dynamics with inherent communication delays are discussed. Section-\ref{sec:simulation} provides three example cases demonstrating significant outcomes of the proposed methodology, and with Section-\ref{sec:conclusions}, the paper is concluded.

\section{PRELIMINARIES AND PROBLEM STATEMENT}\label{sec:prob_formulation}

In this section, we introduce the definitions and notations from graph theory and matrix theory. Then we formulate the main problem to be studied. 

In this work, $R$ denotes the real space. For a real matrix $A$, its transpose and inverse are denoted as $A^T$ and $A^{-1}$, respectively. The symbol $\otimes$ represents the Kronecker product. For matrices $X$ and $P$, the Euclidean norm of $X$ is denoted by $\|X\|$ and the $P$-weighted norm of $X$ is denoted by $\|X\|_P=\sqrt[P]{X^TPX}$. $I_n$ stands for the identity matrix of dimension $n$. Given a matrix $P$, $P>0$ ($P<0$) represents that the matrix is positive definite (or negative definite). Here, we define the column operation $\text{col}(x_1,x_2,\cdots,x_n)$ as $(x^T_1,x^T_2,\cdots,x^T_n)^T$ where $x_1,x_2,\cdots,x_n$ are column vectors.

Consider a multi-agent system of $M$ nonlinear agents. For each agent $i$, the dynamic system is given by:
\begin{equation}\label{model}
\dot{\mathbf{x}}_i=f(\mathbf{x}_i,\mathbf{x}_{-i},t)=\mathbf{F}(x_i)+\mathbf{u}_i,
\end{equation}
where $\mathbf{x}_{i}=(x_{i1},x_{i2},\cdots,x_{in})^{T}$ is the
state vector of the $i$th oscillator, $\mathbf{x}_{-i}$ are the collection of agent $i$'s neighbor's states, the function $\mathbf{F}(\cdot)$ is the corresponding nonlinear vector field, and $\mathbf{u}_i$ is the control input of agent $i$. Here, function $\mathbf{F}(\cdot)$ satisfies the global Lipschitz condition. Therefore there exists positive constant $\beta_i$ such that
\begin{equation*}
\|\mathbf{F}(\mathbf{x}_i)-\mathbf{F}(\mathbf{x}_j)\|\leq \beta_i\|\mathbf{x}_i-\mathbf{x}_j\|.
\end{equation*}
This condition is satisfied if the Jacobians $\frac{\partial F_i}{\partial x_i}$ are uniformly bounded.

There exists a communication network among these agents and the network can be described as an undirected or directed graph $\mathcal{G}=(\mathcal{V,E,C})$. Here $\mathcal{V}=\{1,2,\cdots,M\}$ denotes the node set and $\mathcal{E}\subset \mathcal{V}\times \mathcal{V}$ denotes the edge set. $\mathcal{A}=[a_{ij}]\in R^{M\times M}$ is the adjacency matrix. In this framework, if there exists  a connection between $i$ and $j$ nodes(agents), then $a_{ij}>0$; otherwise, $a_{ij}=0$. We assume there is no self-circle in the graph $\mathcal{G}$, i.e., $a_{ii}=0$. A path is a sequence of connected edges in a graph. If there is a path between any two nodes, the graph is said to be connected. If $\mathcal{A}$ is a symmetric matrix,  $\mathcal{G}$ is called an undirected graph. The set of neighbors of node $i$ is denoted by $\mathcal{N}_i=\{j|(i,j)\subset \mathcal{E}\}$.  The in-degree of agent $i$ is denoted as $\text{deg}_i=\sum^{M-1}_{j=1}a_{ij}$ and the degree matrix is denoted as $\mathcal{D}=\text{diag}(\text{deg}_1,\cdots,\text{deg}_M)$. The Laplacian matrix of $\mathcal{G}$ is described as $\mathcal{L=D-A}$.

In this paper, we consider the multi-agent network systems with switching topologies, where interconnected structures of the network vary with respect to time. For given formulation, denote $\mathcal{P}$ as an index set and $\sigma(t):[0,\infty) \rightarrow \mathcal{P}$ be a switching signal that is defined as a piecewise function. At each time $t$, the graph is represented as $\mathcal{G}_{\sigma(t)}$ and $\{\bar{\mathcal{G}}_p | p \in \mathcal{P}\}$ include all possible graph on the node set $\mathcal{V}$. $\bar{\mathcal{G}}$ is denoted as union of all possible subgraphs $\mathcal{G}_{\sigma(t)}$. Here, each subgraph may be disconnected, where the network becomes jointly connected if the collection $\bar{\mathcal{G}}$ contains a spanning tree. In the light of these, we provide the following definition:

\textbf{Definition-1:} The nonlinear multi-agent network system, given in \eqref{model}, is said to achieve consensus if,
\begin{equation}
\lim_{t\rightarrow \infty}=\|\mathbf{x}_i(t)-\mathbf{x}_j(t)\|=0, \, j=1,\cdots,M,
\end{equation}
is satisfied under switching topology $\mathcal{G}_{\sigma(t)}$, with a distributed control protocol $\mathbf{u}_i=\mu(\mathbf{x}_i,\mathbf{x}_{-i})$. Here $\mathbf{x}_{-i}$'s are the collection of agent $i$'s neighbor's states (i.e., $\mathbf{x}_{-i}=\{\mathbf{x}_{j},j\in \mathcal{N}_i\}$).

In such formulation, the main goal of this study becomes to design a real-time nonlinear receding horizon control based distributed control strategy $\mathbf{u}_i=\mu(\mathbf{x}_i,\mathbf{x}_{-i})$,for each agent $i$, that will achieve consensus, within the given switching network topology and geometric constraints.

\section{DISTRIBUTED NONLINEAR RECEDING HORIZON CONTROL PROTOCOL}\label{sec:dist_control}

In this framework, the following optimization problem is utilized to generate the \emph{local} consensus protocol within the given network, for each specific agent $i$:

\textbf{Problem-1:}
\begin{equation}
\mathbf{u}^*_i(t)=\argmin_{\mathbf{u}_i(t)}J_i(\mathbf{x}_i(t),\mathbf{u}_i(t),\mathbf{x}_{-i}(t))
\end{equation}
\hskip 70pt subject to
\begin{equation*}
\dot{\mathbf{x}}_i(t)=\mathbf{F}(\mathbf{x}_i(t))+\mathbf{u}_i(t),
\end{equation*}
where the performance index is designed as follows:
\begin{equation}\label{cost1}
\begin{split}
J_i=&~\varphi_i+\frac{1}{2}\int_t^{t+T}L_i(\mathbf{x}_i,\mathbf{x}_{-i},\mathbf{u}_i),\\
=&\sum_{j\in \mathcal{N}_i} a_{ij}\|\mathbf{x}_i(t+T)-\mathbf{x}_j(t+T)\|^2_{Q_{iN}} \\ &+\frac{1}{2}\int_t^{t+T}(\sum_{j\in \mathcal{N}_i} a_{ij}\|\mathbf{x}_i(\tau)-\mathbf{x}_j(\tau)\|^2_{Q_{i}}+\|\mathbf{u}_i\|^2_{R_i}){\rm d}\tau,
\end{split}
\end{equation}
Here $Q_{iN}>0$,  $Q_{i}>0$ and $R_{i}>0$ are symmetric matrices, and $T$ defines the horizon. In addition, $\varphi_i$ is used to describe(and define) the terminal cost for each agent.

In this problem formulation, a control scheme is utilized to be able to deal with the nonlinear nature of the topological graph under scrutiny. With this, it is desired to solve the nonlinear optimization problem directly, \emph{in real-time}, without any dependency on iterative processes.

With the construction of the cost function, as given in \eqref{cost1}, the consensus problem is converted into an optimization procedure. For this purpose, we utilize the powerful nature of \emph{real-time nonlinear receding horizon control} algorithm to generate the distributed consensus protocol by minimizing the associated cost function. In this context, each agent only needs to obtain its neighbors' information once via the given network which is more efficient than the centralized control strategy (and the other distributed strategies) that involve multiple information exchanges and predicted trajectories of states. The performance index evaluates the performance from the present time $t$ to the finite future $t+T$, and then is minimized for each time segment $t$ starting from $\mathbf{x}_i(t)$. With this structure, it is possible to convert the present receding horizon control problem into a family of finite horizon optimal control problems on the artificial $\tau$ axis parametrized by time $t$.

According to the first-order necessary conditions of optimality (i.e. for $\delta J_i=0$), a \emph{local} two-point boundary-value problem (TPBVP)  \cite{bryson1975} is formed as follows:
\begin{equation}\label{tpbvp}
\begin{split}
&{\mathbf{\Lambda}_i}^*_{\tau}(\tau,t)=-H^T_{\mathbf{x}_i},\\
&\mathbf{\Lambda}_i^*(T,t)=\varphi_{\mathbf{x}_i}^T[\mathbf{x}^*(T,t)], \\
&{\mathbf{x}_i}^*_{\tau}(\tau,t)=H^T_{\mathbf{\Lambda}_i},\mathbf{x}_i^*(0,t)=\mathbf{x}_i(t),\\
&H_{\mathbf{u}_i}=0.
\end{split}
\end{equation}
where ${\mathbf{\Lambda}_i}$ denotes the costate of each agent $i$ and $H_i$ is the Hamiltonian which is defined as
\begin{equation}\label{eq:tpbvp_ham}
\begin{split}
H_i&= L_i + \mathbf{\Lambda}_{i}^{*T}\dot{\mathbf{x}}_i\\
&=\frac{1}{2}(\sum_{j\in \mathcal{N}_i} a_{ij}\|\mathbf{x}_i(\tau)-\mathbf{x}_j(\tau)\|^2_{Q_{i}}+\|\mathbf{u}_i\|^2_{R_i})+ \mathbf{\Lambda}_{i}^{*T}(\mathbf{F}(x_i)+\mathbf{u}_i).
\end{split}
\end{equation}
Then we have
\begin{equation}\label{tpbvp_l}
\begin{split}
&{\mathbf{\Lambda}_i}^*_{\tau}(\tau,t)=-[\sum_{j\in \mathcal{N}_i} Q_{i}a_{ij}(\mathbf{x}_i(\tau)-\mathbf{x}_j(\tau))+\mathbf{\Lambda}^T_i\mathbf{F}_{\mathbf{x}_i}(\mathbf{x}_i)],\\
&\mathbf{\Lambda}_i^*(T,t)=\sum_{j\in \mathcal{N}_i} Q_{iN}a_{ij}(\mathbf{x}_i(\tau)-\mathbf{x}_j(\tau)).
\end{split}
\end{equation}
In \eqref{tpbvp}-\eqref{tpbvp_l}, $(~~)^*$ denotes a variable in the optimal control problem so as to distinguish it from its correspondence in the original problem. In this notation, $H_{\mathbf{x}_i}$ denotes the partial derivative of $H$ with respect to ${\mathbf{x}_i}$, and so on.

In this methodology, since the state and co-state at $\tau=T$ are determined by the TPBVP in Eq.\eqref{tpbvp} from the state and co-state at $\tau=0$, the TPBVP can be regarded as a nonlinear  algebraic equation with respect to the co-state at $\tau=0$ as
\begin{equation}\label{p}
\mathbf{P}_i(\mathbf{\Lambda}_i(t),\mathbf{x}_i(t),T,t)=\mathbf{\Lambda}_i^*(T,t)-\varphi_{\mathbf{x}_i}^T[\mathbf{x}_i^*(T,t)]=0,
\end{equation}
where $\mathbf{\Lambda}_i(t)$ denotes the co-state at $\tau=0$. The actual local control input for each agent is then given by
\begin{equation}\label{u}
\mathbf{u}_i(t)=\text{arg}\{H_{\mathbf{u}_i}[\mathbf{x}_i(t),\mathbf{\Lambda}_i(t),\mathbf{u}_i(t)]=0\}.
\end{equation}

In this formulation, the optimal control $\mathbf{u}_i(t)$ can be calculated directly from Eq.\eqref{u} based on $\mathbf{x}_i(t)$ and $\mathbf{\Lambda}_i(t)$ information, where the ordinary differential equation of $\lambda(t)$ can be solved numerically from Eq.\eqref{p}, in real-time, without any need of an iterative optimization routine. Since the nonlinear equation $\mathbf{P}_i(\mathbf{\Lambda}_i(t),\mathbf{x}_i(t),T,t)$ has to be satisfied at any time $t$, $\frac{{\rm d}\mathbf{P}_i}{{\rm d}t}=0$ holds along the trajectory of the closed-loop system of the receding horizon control. If $T$ is a smooth function of time $t$, it becomes possible to track the solution of $\mathbf{P}_i(\mathbf{\Lambda}_i(t),\mathbf{x}_i(t),T,t)$ with respect to time. However, numerical errors associated with the solution may accumulate as the integration proceeds in practice, and therefore some correction techniques are required to correct such errors in the solution. To address this problem, a stabilized continuation method \cite{kabamba1987,ohtsuka1994,ohtsuka1997,ohtsuka1998} is used. According to this method, it is possible to rewrite the statement as
\begin{equation}\label{cm}
\frac{{\rm d}\mathbf{P}_i}{{\rm d}t}=A_s\mathbf{P}_i,
\end{equation}
where $A_s$ is a Hurwitz matrix, that helps the solution converge to zero, exponentially.

To evaluate the optimal control (by computing derivative of $\mathbf{\Lambda}_i(t)=\mathbf{\Lambda}_i^*(0,t)$) in real time, we consider the partial differentiation of \eqref{tpbvp} with respect to time $t$,

\begin{equation}
\begin{split}
&\delta \dot{\mathbf{x}}_i=f_{\mathbf{x}_i}\delta \mathbf{x}_i+f_{\mathbf{u}_i}\delta \mathbf{u}_i,\\
&\delta \dot{\mathbf{\Lambda}}_i=-H_{\mathbf{x}_i\mathbf{x}_i}\delta \mathbf{x}_i-H_{\mathbf{x}_i\mathbf{\Lambda}_i}\delta \mathbf{\Lambda}_i-H_{\mathbf{x}_i\mathbf{u}_i}\delta \mathbf{u}_i\\
&0=H_{\mathbf{u}_i\mathbf{x}_i}\delta \mathbf{x}_i+f_{\mathbf{u}_i}^T\delta \mathbf{\Lambda}_i+H_{\mathbf{u}_i\mathbf{u}_i}\delta \mathbf{u}_i.
\end{split}
\end{equation}

Since $\delta \mathbf{u}_i=-H^{-1}_{\mathbf{u}_i\mathbf{u}_i}(H_{\mathbf{u}_i\mathbf{x}_i}\delta \mathbf{x}_i+f_{\mathbf{u}_i}^T\delta \mathbf{\Lambda}_i)$, we have
\begin{align*}
&\delta \dot{\mathbf{x}}_i=(f_{\mathbf{x}_i}-f_{\mathbf{u}_i}H^{-1}_{\mathbf{u}_i\mathbf{u}_i}H_{\mathbf{u}_i\mathbf{u}_i})\delta \mathbf{x}_i-f_{\mathbf{u}_i}H^{-1}_{\mathbf{u}_i\mathbf{u}_i}f_{\mathbf{u}_i}^T\delta \mathbf{\Lambda}_i,\\
&\delta \dot{\mathbf{\Lambda}}_i=-(H_{\mathbf{x}_i\mathbf{x}_i}-H_{\mathbf{x}_i\mathbf{u}_i}H^{-1}_{\mathbf{u}_i\mathbf{u}_i}H_{\mathbf{u}_i\mathbf{x}_i})\delta \mathbf{x}_i-(f_{\mathbf{x}_i}^T-H_{\mathbf{x}_i\mathbf{u}_i}H^{-1}_{\mathbf{u}_i\mathbf{u}_i}f_{\mathbf{u}_i}^T)\delta \mathbf{\Lambda}_i,
\end{align*}  

which leads to the following form of a linear differential equation:
\begin{equation}\label{pd}
\frac {\partial}{\partial \tau}\begin{bmatrix}{\mathbf{x}_i}^*_t-{\mathbf{x}_i}^*_{\tau}\\
{\mathbf{\Lambda}_i}^*_t-{\mathbf{\Lambda}_i}^*_{\tau} \end{bmatrix}=\begin{bmatrix}A_i&-B_i\\
-C_i&-A_i^T \end{bmatrix}\begin{bmatrix}{\mathbf{x}_i}^*_t-{\mathbf{x}_i}^*_{\tau}\\
{\mathbf{\Lambda}_i}^*_t-{\mathbf{\Lambda}_i}^*_{\tau} \end{bmatrix}
\end{equation}
where
$$A_i=f_{\mathbf{x}_i}-f_{\mathbf{u}_i}H^{-1}_{{\mathbf{u}_i}{\mathbf{u}_i}}H_{{\mathbf{u}_i}{\mathbf{x}_i}},$$
$$B_i=f_{\mathbf{u}_i}H^{-1}_{{\mathbf{u}_i}{\mathbf{u}_i}}f_{{\mathbf{u}_i}}^T,$$
$$C_i=H_{{\mathbf{x}_i}{\mathbf{x}_i}}-H_{{\mathbf{x}_i}{\mathbf{u}_i}}H^{-1}_{{\mathbf{u}_i}{\mathbf{u}_i}}H_{{\mathbf{u}_i}{\mathbf{x}_i}}.$$ And the matrix $H_{\mathbf{u}_i\mathbf{u}_i}$ should be non-singular.

In order to reduce the computational cost without resorting to any approximation technique, the backward-sweep method \cite{bryson1975,ohtsuka1997,ohtsuka1998} is implemented. The derivative of the function $\mathbf{P}_i$ with respect to time is given by
\begin{equation}\label{df}
\begin{split}
\frac{{\rm d}\mathbf{P}_i}{{\rm d}t}=&{\mathbf{\Lambda}_i}^*_t(T,t)-\varphi_{{\mathbf{x}_i}{\mathbf{x}_i}}{\mathbf{x}_i}^*_t(T,t)+[{\mathbf{\Lambda}_i}^*_{\tau}(T,t)-\varphi_{{\mathbf{x}_i}{\mathbf{x}_i}}{\mathbf{x}_i}^*_{\tau}(T,t)]\frac{{\rm d}T}{{\rm d}t},
\end{split}
\end{equation}
where $x^*_{\tau}$ and ${\mathbf{\Lambda}_i}^*_{\tau}$ are given by \eqref{tpbvp}.

The relationship between the co-state and other variables is assumed as follows:
\begin{equation}\label{relation}
{\mathbf{\Lambda}_i}^*_t-{\mathbf{\Lambda}_i}^*_\tau=S_i(\tau,t)({\mathbf{x}_i}^*_t-{\mathbf{x}_i}^*_\tau)+c_i(\tau,t).
\end{equation}
which leads to
\begin{equation}\label{sct}
\begin{split}
&S_i(T,t)=\varphi_{\mathbf{x}_i\mathbf{x}_i}\mid_{\tau=T},\\
&c_i(T,t)=(H_{\mathbf{x}_i}^T+\varphi_{\mathbf{x}_i\mathbf{x}_i}f)\mid_{\tau=T}(1+\frac{{\rm d}T}{{\rm d}t})+A_s\mathbf{P}_i.
\end{split}
\end{equation}
According to \eqref{relation} and \eqref{pd}, it becomes possible to form the following differential equations:
\begin{equation}\label{sc}
\begin{split}
&\frac {\partial {S_i}}{\partial \tau}=-A_i^TS_i-S_iA_i+S_iB_iS_i-C_i,\\
&\frac {\partial {c_i}}{\partial \tau}=-(A_i^T-S_iB_i)c_i.
\end{split}
\end{equation}
Based on \eqref{relation}, the differential equation of the co-state to be integrated in real time is obtained as:
\begin{align}\label{dl}
\frac{{\rm d}\mathbf{\Lambda}_i(t)}{{\rm d}t}=-H_{\mathbf{x}_i}^T+c_i(0,t).
\end{align}


The NRHC method for computing the distributed optimal control $\mathbf{u}_i(t)$ is summarized in Algorithm-1, where $t_s >0$ denotes the sampling time.

\begin{algorithm}\label{alg_1}
 \caption{ $\mathbf{u}^*_i(t)=\argmin_{\mathbf{u}_i(t)}J_i(\mathbf{x}_i(t),\mathbf{u}_i(t),\mathbf{x}_{-i}(t))$} \label{}
 (1) Set $t=0$ and initial state $\mathbf{x}_i=\mathbf{x}_i(0)$.\\
 (2) For $t^{\prime}\in[t,t+t_s]$, integrate the defined TPBVP in \eqref{tpbvp_l} forward from $t$ to $t+T$, then integrate \eqref{sc} backward with terminal conditions provided in \eqref{sct} from $t+T$ to $t$. \\
 (3) Integrate the differential equation of $\mathbf{\Lambda}_i(t)$, from $t$ to $t+ t_s$.\\
 (4) At time $t+t_s$, compute $\mathbf{u}^*_i$ by Eq.\eqref{u} with the terminal values of $\mathbf{x}_i(t)$ and $\mathbf{\Lambda}_i(t)$.\\
 (5) Set $t=t+t_s$, return to Step-(2).

\end{algorithm}

\textbf{Lemma-1:} The cost function, defined in Eq. \eqref{cost1}, is strictly convex and guarantees the global minimum.

\emph{Proof}: Since all weighting functions maintain positive definite nature in their structure (such as $Q_{iN}>0$,  $Q_{i}>0$, $R_{i}>0$), from the Karush-Kuhn-Tucker(KKT) conditions \cite{antoniou2007}, the proposed method guarantees the global minima. \qedsymbol

\section{CONVERGENCE AND STABILITY ANALYSIS}\label{sec:stab_analysis}
For the sake of clarity, and without loss of generality, here we define the consensus error as
$$\delta_1(t)=\mathbf{x}_i(t)-\mathbf{x}_1(t)$$
for all $i$ and $\Delta(t)=\text{col}(\delta_1(t),\delta_2(t),\cdots,\delta_M(t))$. The optimal control $\mathbf{U}$ is denoted as $\mathbf{U}=\text{col}(\mathbf{u}_1,\cdots,\mathbf{u}_M)$ for all $i$. 

The cost function can be written as
\begin{equation}
J=\mathbf{\Phi}+\frac{1}{2}\int_t^{t+T}[\Delta^{*T}Q\Delta^*+\mathbf{U}^{*T}R\mathbf{U}^*]{\rm d}\tau,
\end{equation}
where $\mathbf{\Phi}=\sum_{i \in M}\varphi_i$, $Q=\text{col}(Q_1, \cdots, Q_M)$ and $R=\text{col}(R_1, \cdots, R_M)$.

In order to ensure the closed-loop stability of the proposed nonlinear receding horizon control scheme, we first consider the case that terminal cost $\mathbf{\Phi}=0$ and introduce following definitions.

In this regard, we assume the sub-level sets 
\begin{equation*}
\Gamma_r^{\infty}=\{\Delta\in\Gamma^{\infty}:J^*_{\infty}<r^2\}
\end{equation*}
are compact and path connected where $J^*_{\infty}=\int_0^{\infty}[\Delta^{*T}Q\Delta^*+\mathbf{U}^{*T}R\mathbf{U}^*]{\rm d}\tau$ and moreover $\Gamma^{\infty}=\cup_{r\ge 0}\Gamma_r^{\infty}$. We use $r^2$ here to reflect the fact that the cost function is quadratically bounded. And therefore the sub-level set of $\Gamma_r^{T}=\{\Delta\in\Gamma^{\infty}:J^*_{T}<r^2\}$ where $J^*_{T}=\int_t^{t+T}[\Delta^{*T}Q\Delta^*+\mathbf{U}^{*T}R\mathbf{U}^*]{\rm d}\tau$.

\textbf{Lemma-2:} (Dini \cite{jadbabaie2005} ) Let $\{f_n\}$ be a sequence of upper semi-continuous, real-valued functions on a countably compact space $X$, and suppose that for each $x\in X$, the sequence $\{f_n(x)\}$ decreases monotonically to zero. Then the convergence is uniform.

\textbf{Theorem-1:} \cite{jadbabaie2005} Let $r$ be given as $r>0$ and suppose that the terminal cost is equal to zero. For each sampling time $t_s >0$, there exists a horizon window  $T^*<\infty$ such that, for any $T \geqslant T^*$, the receding horizon scheme is asymptotically stabilizing. 

\begin{proof}
By the principle of optimality, we have
\begin{equation*}
J^*_T(\Delta)=\int_{t}^{t+t_s}(\Delta_T^{*T}Q\Delta_T^*+\mathbf{U}_T^{*T} R\mathbf{U}_T^*)d\tau+J^*_{T-t_s}(\Delta_T^*)
\end{equation*}
where $t_s \in [t,t+T]$ is the sampling time and $J^*_{t_s}(\Delta)=\int_{t}^{t+t_s}(\Delta_T^{*T}Q\Delta_T^*+\mathbf{U}_T^{*T} R\mathbf{U}_T^*)d\tau$, so that
\begin{equation*}
\begin{split}
J^*_{T-t_s}(\Delta_T^*)-J^*_{T-t_s}(\Delta)&=J^*_{T}(\Delta)-J^*_{T-t_s}(\Delta)-\int_t^{t+t_s}(\Delta_T^{*T}Q\Delta_T^{*}+\hat{\Theta}_T^{*T} R\hat{\Theta}_T^{*})d\tau\\
&\leq J^*_{t_s}(\Delta)+J^*_{T}(\Delta)-J^*_{T-t_s}(\Delta)
\end{split}
\end{equation*}
Since the terminal cost is equal to zero, it is clear that $T_1<T_2$. This implies that $J^*_{T_1}(\Delta)<J^*_{T_2}(\Delta)$ holds for all $\Delta$ so that
\begin{equation*}
J^*_{T-t_s}(\Delta_T^*)-J^*_{T-t_s}(\Delta) \leq J^*_{t_s}(\Delta)+J^*_{\infty}(\Delta)-J^*_{T-t_s}(\Delta).
\end{equation*}
is satisfied.
If we can show, for example, that there exists a $T^*$ such that $T>T^*$ yields into
\begin{equation*}
J^*_{\infty}(\Delta)-J^*_{T-t_s}(\Delta) \leq \frac{1}{2}J^*_{t_s}(\Delta)
\end{equation*}
for all $\Delta \in \Gamma_r^{\infty}$, stability over any sub-level set of $J^*_{T-t_s}(\cdot)$ that is contained in $\Gamma_r^{\infty}$ will be assured. To that end, define, for $\Delta \in \Gamma_r^{\infty}$

\begin{equation*}
\psi_T(\Delta)=
\begin{cases}
\frac{J^*_{\infty}(\Delta)-J^*_{T-t_s}(\Delta)}{J^*_{t_s}(\Delta)}, & \quad \Delta \neq 0 \\
\lim \sup_{x\rightarrow 0} \psi_T(\Delta), & \quad \Delta=0
\end{cases}
\end{equation*}
where $\psi_T(\cdot)$ is upper semi-continuous on $\Gamma_r^{\infty}$. It is clear that $\psi_T(\cdot)$ is a monotonically decreasing family of upper semi-continuous functions defined over the compact set $\Gamma_r^{\infty}$. Thus, by Dini's theorem (as stated in \cite{jadbabaie2005}), there exists a $T^*<\infty$ such that $\psi_T(\Delta)<\frac{1}{2}$ for all $\Delta \in \Gamma_r^{\infty}$ and all $T\ge T^*$. Here, for each $r_1>0$ we have $\Gamma_{r_1}^{T-t_s}\subset \Gamma_r^{\infty}$ satisfied, leading to

\begin{equation*}
J^*_{T-t_s}(\Delta_T^*)-J^*_{T-t_s}(\Delta) \leq -\frac{1}{2}J^*_{t_s}(\Delta)
\end{equation*}
for all $\Delta \in \Gamma_{r_1}^{T-t_s}$.

\end{proof}

Next, we present the closed-loop stability of the proposed nonlinear receding horizon control scheme with locally quadratic terminal cost, i.e. $\mathbf{\Phi}=\sum_{i \in M}\varphi_i$.

\textbf{Theorem-2:} (\cite{jadbabaie2005}) Let $r$ be given as $r>0$ and suppose that the terminal cost is non-negative and locally quadratically bounded. For each sampling time $t_s>0$, there exists a horizon window  $T^*<\infty$ such that, for any $T \geqslant T^*$, the receding horizon scheme is asymptotically stabilizing. 

\begin{proof}
The theorem utilizes the results of \textbf{Theorem-1}, for the proof. Here, we define $J^*_{T,0}(\cdot)$ to denote the cost function with zero terminal cost and $J^*_{T,1}(\cdot)$ to denote the cost function with locally quadratic terminal cost. It is clear to show that
\begin{equation*}
J^*_{T,0}(\Delta)\leq J^*_{T}(\Delta)\leq J^*_{T,1}(\Delta),
\end{equation*} 
and then 
\begin{equation*}
|J^*_{T}(\Delta)- J^*_{\infty}(\Delta)|\leq \max \{J^*_{\infty}(\Delta)-J^*_{T,0}(\Delta),J^*_{T,1}(\Delta)-J^*_{\infty}(\Delta)\},
\end{equation*} 
for all $\Delta \in \Gamma_r^{\infty}$ so that $J^*_{T}(\cdot)$ also converge uniformly to $J^*_{\infty}(\cdot)$ with respect to any locally quadratic positive definite terminal cost.
\end{proof}

\textbf{Corollary-1:} Consider the nonlinear multi-agent system given in Eq.\eqref{model} and assume that the switching interconnected graph $\bar{\mathcal{G}}$, is jointly connected. For the given distributed control protocol in \textbf{Problem-1}, based on \textbf{Theorem-2}, there exists a large enough value of horizon $T$ which guarantees the consensus error $\Delta$ to remain asymptotically stable to achieve consensus in the multi-agent system.

Although the back-ward sweep algorithm is executable whenever the system is stable or not, with this result, when the optimization horizon is chosen to be sufficiently long, the non-increasing monotonicity of the cost function becomes a sufficient condition for the stability. \emph{Therefore, the stability of the multi-agent nonlinear system under jointly connected switching structure by distributed nonlinear receding horizon control method is also ensured.} \openbox

\section{EXTENSION TO COMMUNICATION TIME-DELAY CASE}\label{sec:time-delay}

In this section, we consider a multi-agent network system consisted of $M$ nonlinear agents with inherent time delay values $t_d>0$, which define associated communication delay characteristics. The underlying dynamics of each agent is given in Eq. \eqref{model}.

At this point of the study, the overall aim is to design the distributed control strategy $\mathbf{u}_i(t)=\mu( \mathbf{x}_i(t),\mathbf{x}_{-i}(t-t_d))$ to achieve consensus within the given network topology in presence of inherent communication delays using the above mentioned real-time nonlinear receding horizon control methodology. 

For each agent $i$, the following optimization problem is utilized to generate the consensus protocol locally, within the network:

\textbf{Problem-2:}
\begin{equation}
\mathbf{u}^*_i(t)=\text{argmin}J_i(\mathbf{x}_i(t),\mathbf{u}_i(t),\mathbf{x}_{-i}(t-t_d))
\end{equation}
\hskip 70pt subject to
\begin{equation*}
\dot{\mathbf{x}}_i(t)=\mathbf{F}(\mathbf{x}_i(t))+\mathbf{u}_i(t),
\end{equation*}
The corresponding performance index is designed as follows:
\begin{equation}\label{cost}
\begin{split}
J_i=&~\varphi_i+\frac{1}{2}\int_t^{t+T}L_i(\mathbf{x}_i(t),\mathbf{x}_{-i}(t-t_d),\mathbf{u}_i(t)),\\
=&\sum_{j\in \mathcal{N}_i} a_{ij}\|\mathbf{x}_i(t+T)-\mathbf{x}_j(t+T)\|^2_{Q_{iN}}\\ &+\frac{1}{2}\int_t^{t+T}(\sum_{j\in \mathcal{N}_i} a_{ij}\|\mathbf{x}_i(\tau)-\mathbf{x}_j(\tau-t_d)\|^2_{Q_{i}}+\|\mathbf{u}_i\|^2_{R_i}){\rm d}\tau.
\end{split}
\end{equation}

We utilize the same framework and algorithm as in Section-\ref{sec:dist_control} to solve the nonlinear consensus problem with communication time-delay, directly (i.e. in real-time). 

\section{NUMERICAL EXAMPLE AND SIMULATION RESULTS}\label{sec:simulation}
In this section, we demonstrate the validity and feasibility of proposed scheme on several multi-agent nonlinear \emph{chaotic} systems.

\begin{figure}
 \centering
 \vspace*{-3.cm} 
    \includegraphics[width=12cm]{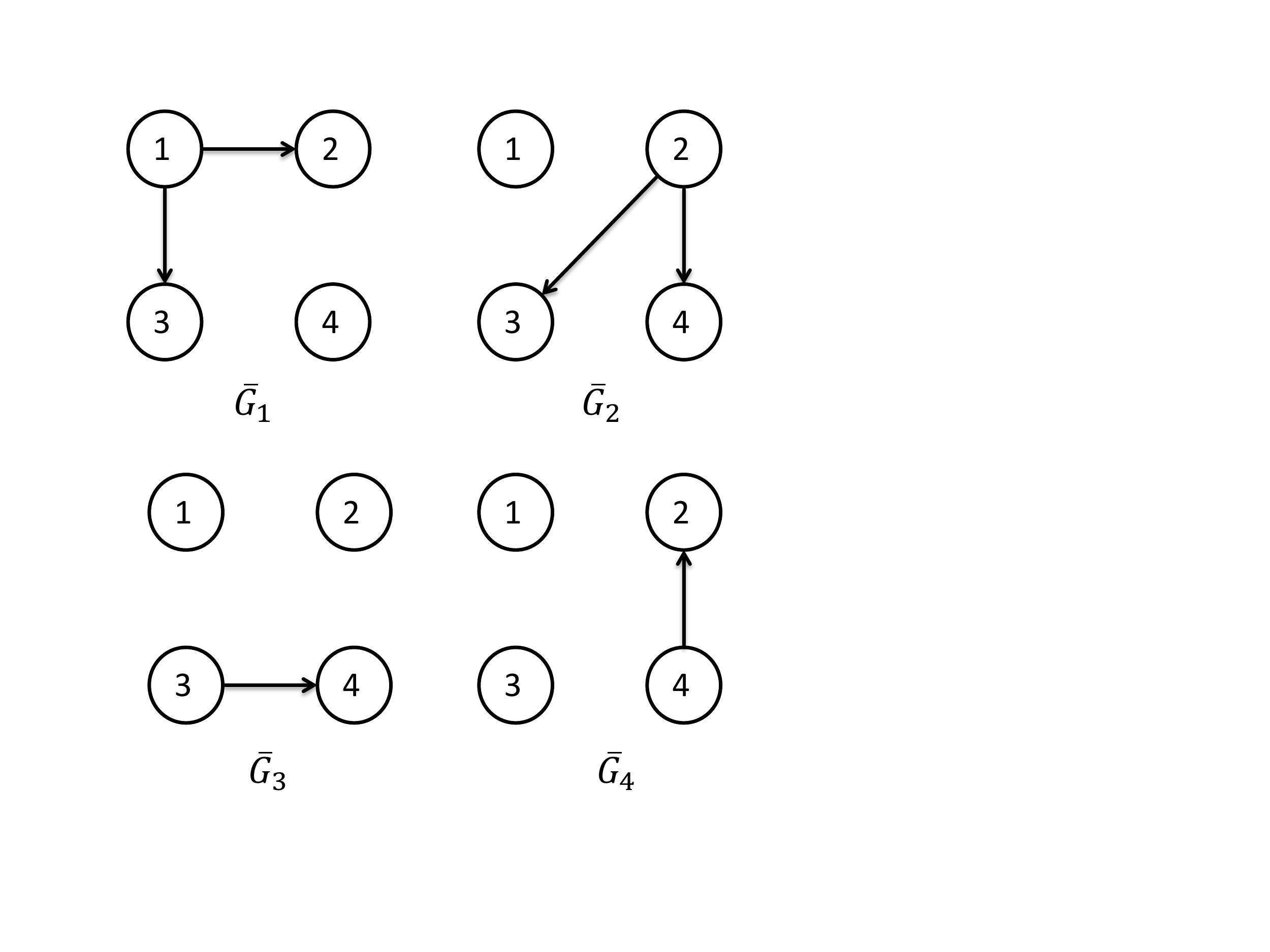}
    \vspace*{-1.1cm}
    \caption{ \label{t1} Four possible directed communication topologies of the multi-agent Lorenz chaotic systems with $4$ agents.}
 \end{figure}

\textbf{Example-1:} First, consider a multi-agent system with $4$ agents, where each agent is modeled as the Lorenz \emph{chaotic} system \cite{lorenz1963}:
\begin{equation}\label{lorenz}
\begin{cases}
\dot{x}_{i1}=10(x_{i2}-x_{i1}), \\
\dot{x}_{i2}=28x_{i1}-x_{i1}x_{i3}-x_{i2},\\
\dot{x}_{i3}=x_{i1}x_{i2}-\frac{8}{3}x_{i3},
\end{cases}
\end{equation}
Here $x_i=(x_{i1},x_{i2},x_{i3})^{T}$ are the states of the $i$-th agent. The switching sub-graph associated with $\sigma(t)$ are  $\bar{\mathcal{G}}_1,\bar{\mathcal{G}}_2,\bar{\mathcal{G}}_3,\bar{\mathcal{G}}_4$, as shown in Fig. \ref{t1}.

Here, the initial states are given by
\begin{equation}\label{initial1}
\begin{bmatrix}x_{i1}(0) \\ x_{i2}(0)\\x_{i3}(0)\end{bmatrix}=\begin{bmatrix}-1&2&-10&9\\ 10& -1& 20& -10\\2& 5& 8& -2\end{bmatrix},
\end{equation}

It can be seen in Fig. \ref{t1} that the switching directed graph $\bar{\mathcal{G}}$  of Lorenz network is jointly connected. The weighting matrices in the cost function are designed as $Q_{iN}=Q_{i}=R_{i}=\text{diag}(1,1,1)$ for all agents. The stable matrix is designed as $A_s=-50I$.

The horizon $T$ in the performance index is given by
\begin{equation}\label{h_t}
T(t)=T_f(1-e^{-\alpha t}),
\end{equation}
where $T_f=1$ and $\alpha=0.01$.

\begin{figure}
 \centering
 \vspace*{-1.1cm}
    \includegraphics[width=13.5cm]{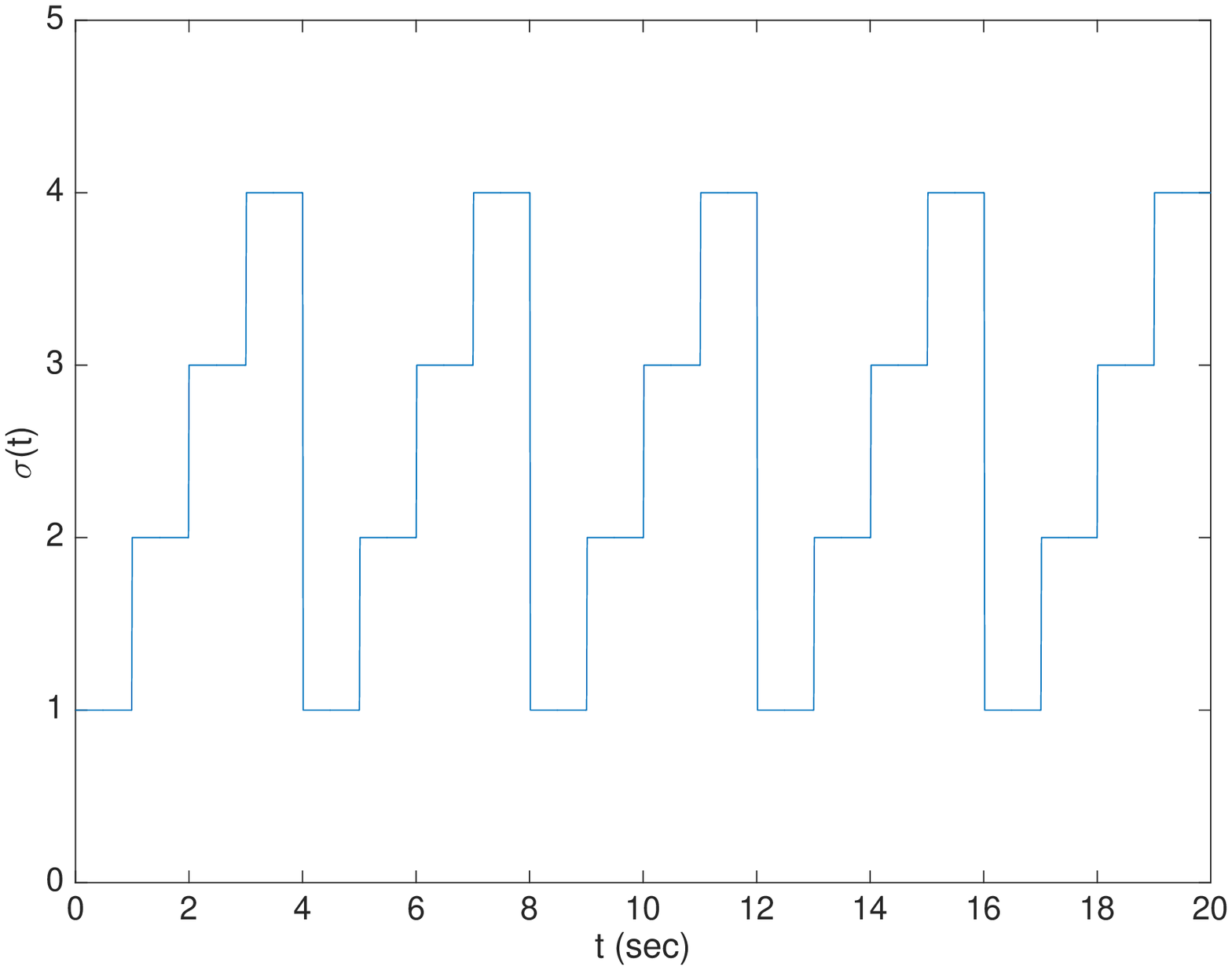}
    \vspace*{-0.5cm}
    \caption{\label{fig_sig1} The switching signal $\sigma(t)$.}
 \end{figure}

\begin{figure}
 \centering
 \vspace*{-1.1cm}
    \includegraphics[width=13.5cm]{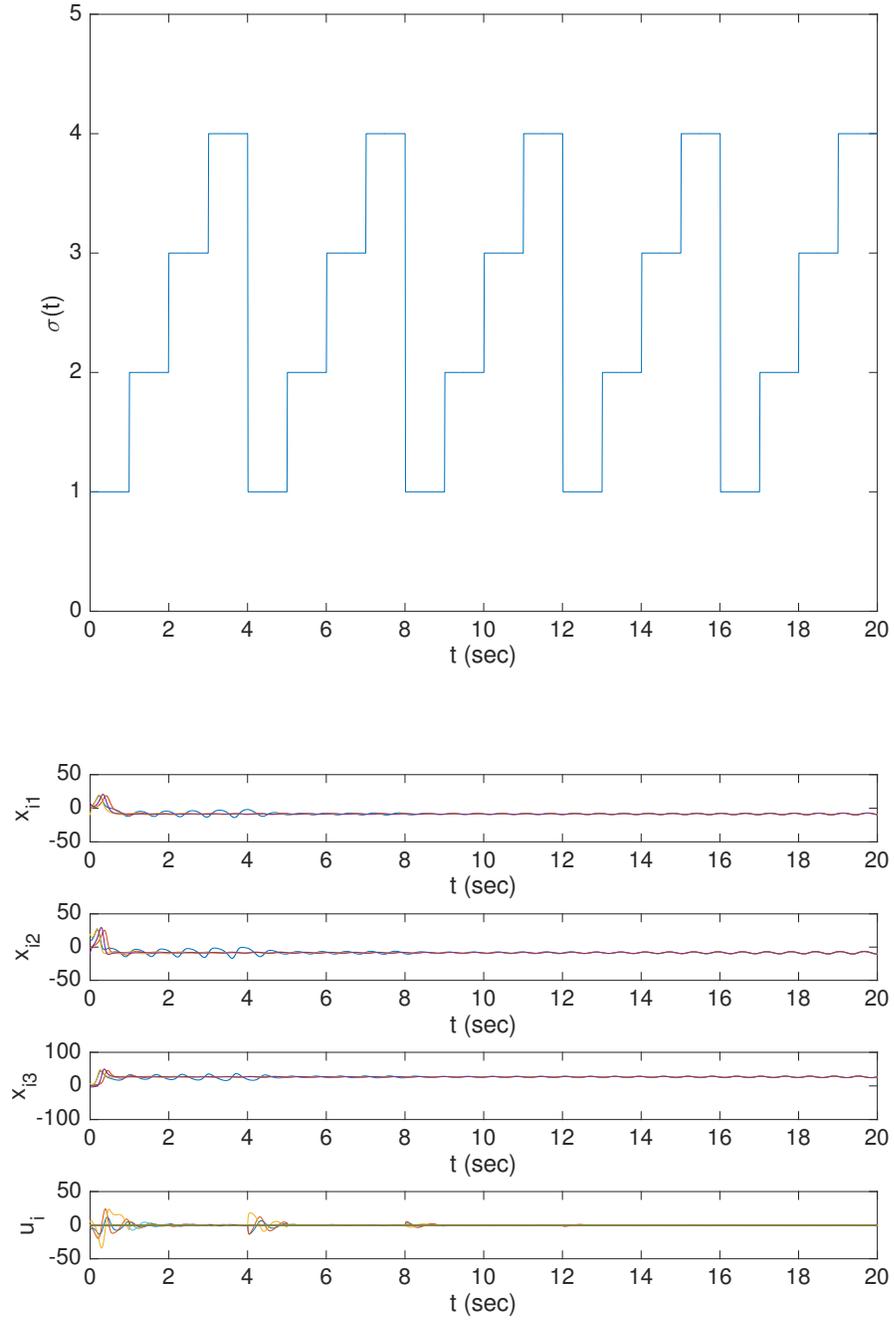}
    \vspace*{-0.5cm}
    \caption{\label{fig_sys1} The trajectories of all agents $x_i(t) (i=1,\cdots,4)$ of Lorenz chaotic system and all control protocol $u_i$ generated by distributed NRHC.}
 \end{figure}

The simulation is implemented in MATLAB, where the sampling time $t_s$ is $0.01s$ and the time step on the artificial $\tau$ axis is $0.005s$. Defined switching signal $\sigma(t)$ is as given in Fig. \ref{fig_sig1}. Fig. \ref{fig_sys1} depicts the trajectories of this multi-agent Lorenz system with initial conditions defined in Eq. \eqref{initial1}. Under the given switching topology, it is possible to observe that the proposed distributed real-time nonlinear receding horizon control strategy in Algorithm-1 results in a consensus, clearly demonstrating the effectiveness of the algorithm. Here, the horizon length is kept sufficiently long to ensure the stability.

\begin{figure}
 \centering
 \vspace*{-1.cm}
    \includegraphics[width=12cm]{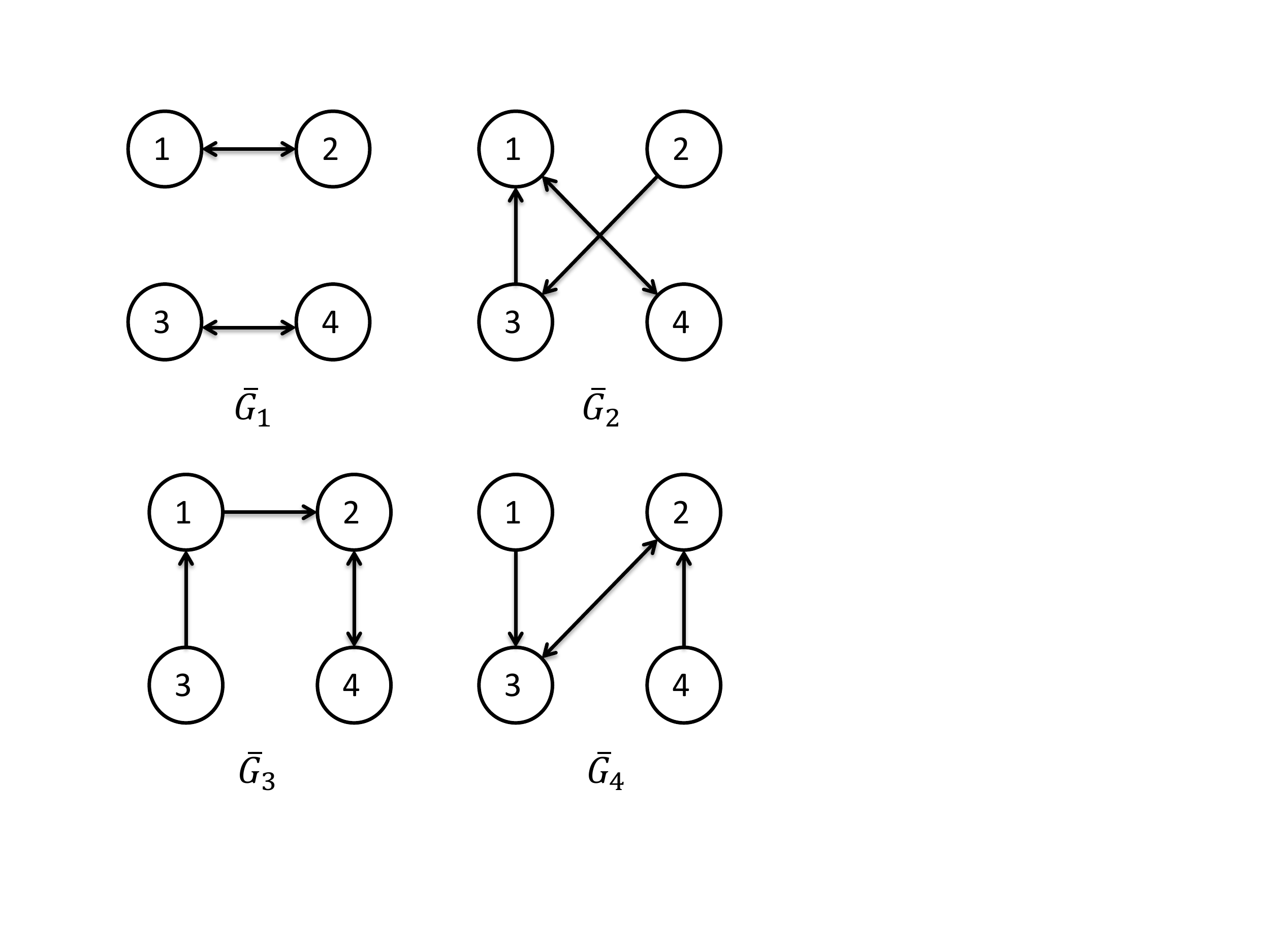}
    \vspace*{-1.1cm}
    \caption{ \label{t2} Four possible directed communication topologies of the multi-agent Lorenz chaotic systems with $4$ agents.}
 \end{figure}

\textbf{Example-2:} Next, consider the same multi-agent system under different switching topologies (as shown in Fig. \ref{t2}), where each agent is also modeled as the Lorenz chaotic system.

It can be seen in Fig. \ref{t2} that the graph $\bar{\mathcal{G}}$  is jointly connected. The weighting matrices in the cost function are designed as $Q_{iN}=Q_{i}=R_{i}=\text{diag}(1,1,1)$ for all agents.The stable matrix is designed as $A_s=-50I$. The horizon $T$ in the performance index is given in \eqref{h_t} with $T_f=1$ and $\alpha=0.01$.

Different from the previous example (where \emph{fixed} switching signal was given), in this example, the switching signal is generated \emph{automatically}. In the Step-3 of Algorithm-1, beside the optimal controller, an optimal switching topology is also computed by comparing the minimum of the cost function among all possible topologies. The new algorithm (including the switching logic) is presented in Algorithm-2. 

\begin{algorithm}\label{alg2}
 \caption{ $\mathbf{u}^*_i(t)=\argmin_{\mathbf{u}_i(t)}J_i(\mathbf{x}_i(t),\mathbf{u}_i(t),\mathbf{x}_{-i}(t))$} \label{}
 (1) Set $t=0$, initial state $\mathbf{x}_i=\mathbf{x}_i(0)$ and initial switching signal $\sigma=\sigma(0)$.\\
 (2) For $t^{\prime}\in[t,t+t_s]$, integrate the defined TPBVP in \eqref{tpbvp_l} forward from $t$ to $t+T$, then integrate \eqref{sc} backward with terminal conditions provided in \eqref{sct} from $t+T$ to $t$. \\
 (3) Integrate the differential equation of $\mathbf{\Lambda}_i(t)$, from $t$ to $t+ t_s$, then calculate the cost function among all possible topologies.\\
 (4) At time $t+t_s$,  select the topology associated with the minimum of the cost function as current topology, then compute $\mathbf{u}^*_i$ by Eq.\eqref{u} with the terminal values of $\mathbf{x}_i(t)$ and $\mathbf{\Lambda}_i(t)$.\\
 (5) Set $t=t+t_s$, return to Step-(2).
\end{algorithm}

\begin{figure}
 \centering
 \vspace*{-.1cm}
    \includegraphics[width=13.5cm]{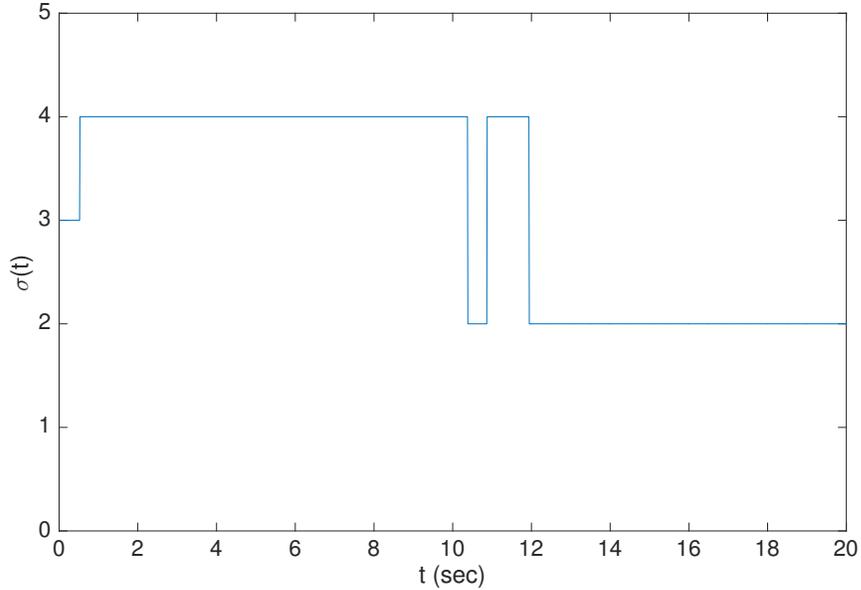}
    \vspace*{-0.5cm}
    \caption{\label{fig_sig2} The switching signal $\sigma(t)$.}
 \end{figure}
 
 \begin{figure}
 \centering
 \vspace*{-.1cm}
    \includegraphics[width=13.1cm]{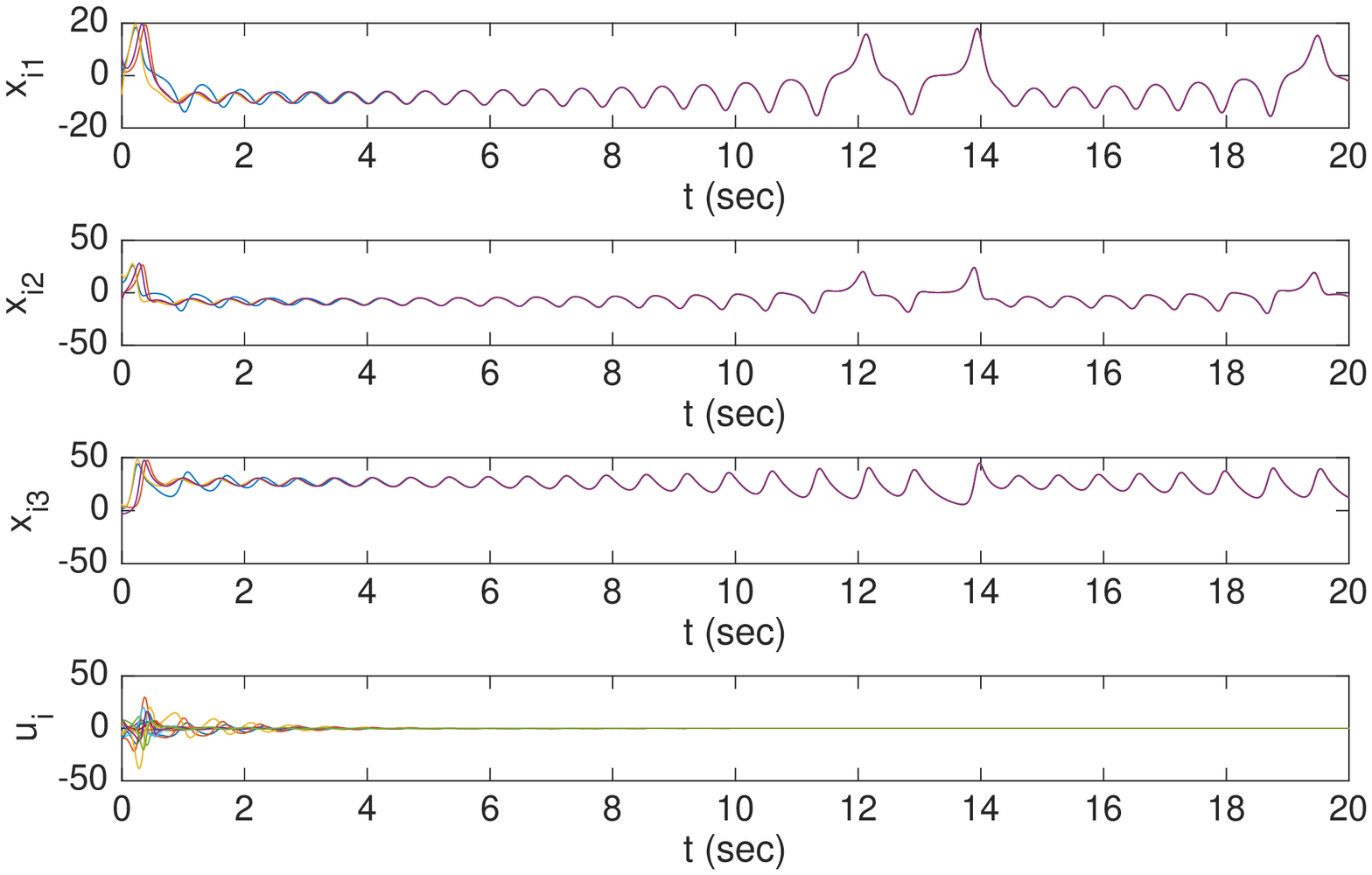}
    \vspace*{-0.5cm}
    \caption{\label{fig_sys2} The trajectories of all agents $x_i(t) (i=1,\cdots,4)$ of Lorenz chaotic system and all control protocol $u_i$ generated by distributed NRHC.}
 \end{figure}

The simulation is implemented in MATLAB, where the sampling time $t_d$
is $0.01s$ and the time step on the artificial $\tau$ axis is $0.005s$. Fig. \ref{fig_sys2} depicts the trajectories of this multi-agent systems in \eqref{lorenz} with initial conditions defined in Eq. \eqref{initial1} and the corresponding control strategies under the optimal switching signal $\sigma(t)$ as Fig. \ref{fig_sig2}, where the consensus is achieved through the suggested distributed real-time nonlinear receding horizon control method.

\textbf{Example-3:} In this scenario, let's consider the same (as in Example-2) Lorenz system with inherent communication delay(s). For this case, the adjacency matrix $\mathcal{A}$ of $\mathcal{G}$ is given as
\begin{equation}
\begin{bmatrix}0&0&1&0 \\ 1&0&1&0\\0&1&0&1\\1&0&0&0\end{bmatrix}
\end{equation}
The weighting matrices in the cost function are designed as $Q_{i0}=\text{diag}(10,10,10)$ and $Q_{iN}=Q_{i}=R_{i}=\text{diag}(1,1,1)$ for all agents. The stable matrix is designed as $A_s=-50I$. The horizon $T$ in the performance index is given in \eqref{h_t} with $T_f=1$ and $\alpha=0.01$.

 \begin{figure}
 \centering
 \vspace*{-.1cm}
    \includegraphics[width=12.3cm]{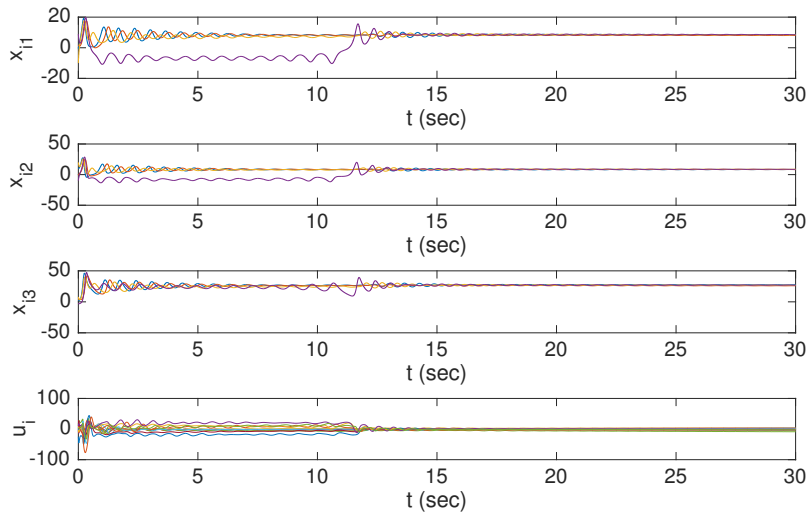}
    \vspace*{-0.5cm}
    \caption{\label{fig_sys3} The trajectories of all agents $x_i(t) (i=1,\cdots,4)$ of Lorenz chaotic system and all control protocol $u_i$ generated by distributed NRHC with communication time-delay $t_d=0.2$ sec.}
 \end{figure}
 
\begin{figure}
 \centering
 \vspace*{-.1cm}
    \includegraphics[width=12.3cm]{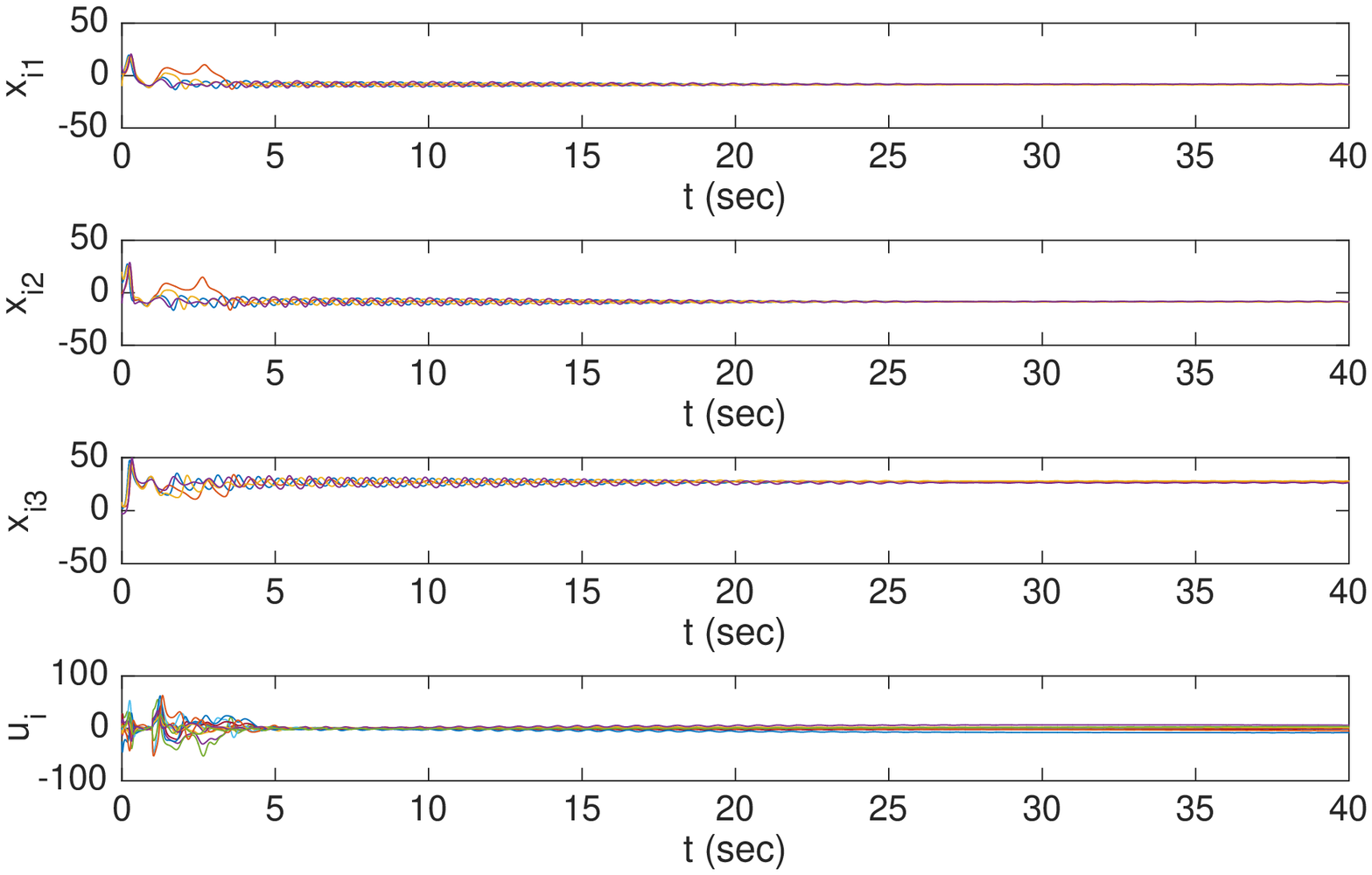}
    \vspace*{-0.5cm}
    \caption{\label{fig_sys4} The trajectories of all agents $x_i(t) (i=1,\cdots,4)$ of Lorenz chaotic system and all control protocol $u_i$ generated by distributed NRHC with communication time-delay $t_d=1$ sec.}
 \end{figure}

The simulation is implemented in MATLAB, where the sampling time $t_s$ is $0.01s$ and the time step on the artificial $\tau$ axis is $0.005s$. Fig. \ref{fig_sys3} depicts the trajectories of this multi-agent Lorenz system with initial condition in Eq. \eqref{initial1} under communication time-delay $t_d=0.2$ sec and the corresponding control strategies where the consensus is reached by using the distributed real-time nonlinear receding horizon control method. Fig. \ref{fig_sys4} depicts the trajectories of this system with communication time-delay $t_d=1$ sec.

As it can be seen easily from above example, proposed real-time nonlinear receding horizon control methodology is working remarkable within switching topology and communication time-delay. All the agents in the systems are able to reach consensus. Here, again, the horizon length is kept sufficiently long to ensure the stability.

\section{CONCLUSIONS }\label{sec:conclusions}

In this paper, we investigated the multi-agent consensus problem of nonlinear systems under switching topologies and embedded communication time-delay by using distributed real-time nonlinear receding horizon control methodology. Different from the previous works, we solved the nonlinear optimal consensus problem directly, without any need or the utilization of linearization techniques and/or iterative procedures. Based on the stabilized continuation method, the backward sweep algorithm is implemented to minimize the consensus error among the agents and the local control strategy is integrated in real time. We provided stability guarantees of the systems if the horizon length is kept sufficiently long. Several benchmark examples with different topologies demonstrates the applicability and significant outcomes of proposed scheme on nonlinear chaotic systems. 

\section*{References}

%
%
%
%

\end{document}